
\magnification=\magstep1
\baselineskip =5mm
\lineskiplimit =1.0mm
\lineskip =1.0mm

\long\def\comment#1{}

\long\def\blankout #1\eb{}
\def\noblankout{\def\blankout{}\def\eb{}}

\let\properlbrack=\lbrack
\let\properrbrack=\rbrack
\def\ordcomma{,}
\def\ordcolon{:}
\def\ordsemicolon{;}
\def\ordleftparen{(}
\def\ordrightparen{)}
\def\ordleftbrack{\properlbrack}
\def\ordrightbrack{\properrbrack}
\def\rmcomma{\ifmmode ,\else \/{\rm ,}\fi}
\def\rmcolon{\ifmmode :\else \/{\rm :}\fi}
\def\rmsemicolon{\ifmmode ;\else \/{\rm ;}\fi}
\def\rmleftparen{\ifmmode (\else \/{\rm (}\fi}
\def\rmrightparen{\ifmmode )\else \/{\rm )}\fi}
\def\rmleftbrack{\ifmmode \properlbrack\else \/{\rm \properlbrack}\fi}
\def\rmrightbrack{\ifmmode \properrbrack\else \/{\rm \properrbrack}\fi}
\catcode`,=\active 
\catcode`:=\active 
\catcode`;=\active 
\catcode`(=\active 
\catcode`)=\active 
\catcode`[=\active 
\catcode`]=\active 
\let,=\ordcomma
\let:=\ordcolon
\let;=\ordsemicolon
\let(=\ordleftparen
\let)=\ordrightparen
\let[=\ordleftbrack
\let]=\ordrightbrack
\let\lbrack=\ordleftbrack
\let\rbrack=\ordrightbrack
\def\rmpunctuation{
\let,=\rmcomma
\let:=\rmcolon
\let;=\rmsemicolon
\let(=\rmleftparen
\let)=\rmrightparen
\let[=\rmleftbrack
\let]=\rmrightbrack
\let\lbrack=\rmleftbrack
\let\rbrack=\rmrightbrack}

\def\writemonth#1{\ifcase#1
\or January\or February\or March\or April\or May\or June\or July%
\or August\or September\or October\or November\or December\fi}

\newcount\mins
\newcount\minmodhour
\newcount\hour
\newcount\hourinmin
\newcount\ampm
\newcount\ampminhour
\newcount\hourmodampm
\def\writetime#1{%
\mins=#1%
\hour=\mins \divide\hour by 60
\hourinmin=\hour \multiply\hourinmin by -60
\minmodhour=\mins \advance\minmodhour by \hourinmin
\ampm=\hour \divide\ampm by 12
\ampminhour=\ampm \multiply\ampminhour by -12
\hourmodampm=\hour \advance\hourmodampm by \ampminhour
\ifnum\hourmodampm=0 12\else \number\hourmodampm\fi
:\ifnum\minmodhour<10 0\number\minmodhour\else \number\minmodhour\fi
\ifodd\ampm p.m.\else a.m.\fi
}

\font\tenrm=cmr10
\font\smallcaps=cmcsc10
\font\eightrm=cmr8
\font\ninerm=cmr9
\font\sixrm=cmr6
\font\eightbf=cmbx8
\font\sixbf=cmbx6
\font\eightit=cmti8
\font\eightsl=cmsl8
\font\eighti=cmmi8
\font\eightsy=cmsy8
\font\eightex=cmex10 at 8pt
\font\sixi=cmmi6
\font\sixsy=cmsy6
\font\ninesy=cmsy9
\font\seventeenrm=cmr17
 1
\font\twelverm=cmr10 scaled \magstep2
\font\seventeeni=cmmi10 scaled \magstep3
\font\twelvei=cmmi10 scaled \magstep2
\font\seventeensy=cmsy10 scaled \magstep3
\font\twelvesy=cmsy10 at 12pt
\font\seventeenex=cmex10 scaled \magstep3
\font\seventeenbf=cmbx10 scaled \magstep3
\catcode`@=11
\def\eightbig#1{{\hbox{$\textfont0=\ninerm\textfont2=\ninesy
\left#1\vbox to6.5pt{}\right.\n@space$}}}
\catcode`@=12
\def\eightpoint{\eightrm \normalbaselineskip=4.5 mm%
\textfont0=\eightrm \scriptfont0=\sixrm \scriptscriptfont0=\fiverm%
\def\rm{\fam0 \eightrm}%
\textfont1=\eighti \scriptfont1=\sixi \scriptscriptfont1=\fivei%
\def\mit{\fam1 } \def\oldstyle{\fam1 \eighti}%
\textfont2=\eightsy \scriptfont2=\sixsy \scriptscriptfont2=\fivesy%
\def\cal{\fam2 }%
\textfont3=\eightex \scriptfont3=\eightex \scriptscriptfont3=\eightex%
\def\bf{\fam\bffam\eightbf} \textfont\bffam\eightbf
\scriptfont\bffam=\sixbf \scriptscriptfont\bffam=\fivebf
\def\it{\fam\itfam\eightit} \textfont\itfam\eightit
\def\sl{\fam\slfam\eightsl} \textfont\slfam\eightsl
\let\big=\eightbig \normalbaselines\rm
\def\caps##1{\nottencaps{##1}}
}
\def\seventeenpoint{\seventeenrm \baselineskip=5.5mm%
\textfont0=\seventeenrm \scriptfont0=\twelverm \scriptscriptfont0=\sevenrm%
\def\rm{\fam0 \seventeenrm}%
\textfont1=\seventeeni \scriptfont1=\twelvei \scriptscriptfont1=\seveni%
\def\mit{\fam1 } \def\oldstyle{\fam1 \seventeeni}%
\textfont2=\seventeensy \scriptfont2=\twelvesy \scriptscriptfont2=\sevensy%
\def\cal{\fam2 }%
\textfont3=\seventeenex \scriptfont3=\seventeenex%
\scriptscriptfont3=\seventeenex%
\def\bf{\fam\bffam\seventeenbf} \textfont\bffam\seventeenbf
}

\def\setheadline #1\\ #2 \par{\headline={\ifnum\pageno=1 
\hfil
\else \eightpoint \noindent
\ifodd\pageno \hfil \caps{#2}\hfil \else
\hfil \caps{#1}\hfil \fi\fi}}

\footline={\ifnum\pageno=1 \hfil
\else \tenrm \hfil \folio \hfil \fi}

\def\beginsection{} 
\def\datedversion{\footline={\ifnum\pageno=1 \fiverm \hfil
Typeset using plain-\TeX\ on
\writemonth\month\ \number\day, \number\year\ at \writetime{\time}\hfil 
\else \tenrm \hfil \folio \hfil \fi}
\def\tempsetheadline##1{\headline={\ifnum\pageno=1 
\hfil
\else \eightpoint \noindent
\writemonth\month\ \number\day, \number\year,
\writetime{\time}\hfil ##1\fi}}
\def\firstbeginsection##1\par{\bigskip\vskip\parskip
\message{##1}\centerline{\caps{##1}}\nobreak\smallskip\noindent
\tempsetheadline{##1}} \def\beginsection##1\par{\vskip0pt
plus.3\vsize\penalty-250 \vskip0pt plus-.3\vsize\bigskip\vskip\parskip
\message{##1}\centerline{\caps{##1}}\nobreak\smallskip\noindent
\tempsetheadline{##1}}}

\def\finalversion{\footline={\ifnum\pageno=1 \eightrm \hfil 
This paper is in final form.\hfil 
\else \tenrm \hfil \folio \hfil \fi}}

\def\preliminaryversion{\footline={\ifnum\pageno=1 \eightrm \hfil 
Preliminary Version.\hfil 
\else \tenrm \hfil \folio \hfil \fi}}

\def\Head #1: {\medskip\noindent{\it #1}:\enspace}
\def\Proof: {\Head Proof: }
\def\Proofof #1: {\Head Proof of #1: }
\def\endproof{\nobreak\hfill$\sqr$\bigskip\goodbreak}

\def\Abstract\par#1\par{\centerline{\vtop{
\eightpoint
\abovedisplayskip=6pt plus 3pt minus 3pt
\belowdisplayskip=6pt plus 3pt minus 3pt
\moreabstract\parindent=0 true in%
\caps{Abstract}: \ \ #1}}
\abovedisplayskip=12pt plus 3pt minus 9pt
\belowdisplayskip=12pt plus 3pt minus 9pt
\vskip 0.4 true in}
\def\moreabstract{%
\par \hsize = 5 true in \hangindent=0 true in \parindent=0.5 true in}

\def\caps#1{\smallcaps #1}
\def\nottencaps#1{\uppercase{#1}}

\def\firstbeginsection#1\par{\bigskip\vskip\parskip
\message{#1}\centerline{\caps{#1}}\nobreak\smallskip\noindent}

\def\beginsection#1\par{\vskip0pt plus.3\vsize\penalty-250
\vskip0pt plus-.3\vsize\bigskip\vskip\parskip
\message{#1}\centerline{\caps{#1}}\nobreak\smallskip\noindent}

\def\proclaim#1. #2\par{
\medbreak
\noindent{\caps{#1}.\enspace}{\it\rmpunctuation#2\par}
\ifdim\lastskip<\medskipamount \removelastskip
\penalty55\medskip\fi}

\def\Definition: #1\par{
\Head Definition: #1\par
\ifdim\lastskip<\medskipamount \removelastskip
\penalty55\medskip\fi}

\def\Problem #1: #2\par{
\Head Problem #1: #2\par
\ifdim\lastskip<\medskipamount \removelastskip
\penalty55\medskip\fi}

\def\sqr{\vcenter {\hrule height.3mm
\hbox {\vrule width.3mm height 2mm \kern2mm
\vrule width.3mm } \hrule height.3mm }}

\def\references#1{{
\frenchspacing
\eightpoint
\rmpunctuation
\halign{\bf##\hfil & \quad\vtop{\hsize=5.5 true
in\parindent=0pt\hangindent=3mm \strut\rm##\strut\smallskip}\cr#1}}}

\def\ref[#1]{{\bf [#1]}}

\catcode`@=11 
\def\vfootnote#1{\insert\footins\bgroup
\eightpoint
\interlinepenalty=\interfootnotelinepenalty
\splittopskip=\ht\strutbox
\splitmaxdepth=\dp\strutbox \floatingpenalty=20000
\leftskip=0pt \rightskip=0pt \spaceskip=0pt \xspaceskip=0pt
\textindent{#1}\footstrut\futurelet\next\fo@t}

\def\footremark{\insert\footins\bgroup
\eightpoint\it\rmpunctuation
\interlinepenalty=\interfootnotelinepenalty
\splittopskip=\ht\strutbox
\splitmaxdepth=\dp\strutbox \floatingpenalty=20000
\leftskip=0pt \rightskip=0pt \spaceskip=0pt \xspaceskip=0pt
\noindent\footstrut\futurelet\next\fo@t}
\catcode`@=12

\def\verses #1\par
{{\eightpoint\it\centerline{\hbox{\vbox{\halign{##\hfil\cr#1}}}}}}

\def\Bbb{\bf}
\def\E{{\Bbb E}}
\def\R{{\Bbb R}}
\def\Z{{\Bbb Z}}
\def\N{{\Bbb N}}
\def\C{{\Bbb C}}
\font\specialeightrm=cmr10 at 8pt
\def\R{\hbox{\rm I\kern-2pt R}}
\def\Z{\hbox{\rm Z\kern-3pt Z}}
\def\N{\hbox{\rm I\kern-2pt I\kern-3.1pt N}}
\def\C{\hbox{\rm \kern0.7pt\raise0.8pt\hbox{\specialeightrm I}\kern-4.2pt C}} 
\def\E{\hbox{\rm I\kern-2pt E}}

\def\list#1,#2{#1_1$, $#1_2,\ldots,$\ $#1_{#2}}

\def\normo#1{{\left\| #1 \right\|}}
\def\snormo#1{{\mathopen\| #1 \mathclose\|}}

\def\Bignormo#1{{\Bigl\| #1 \Bigr\|}}

\def\modo#1{{\left| #1 \right|}}


\def\and{\mathop{\hbox{\ and\ }}}


\noblankout

\setheadline Sums of i.i.d.r.v.\\
             Montgomery-Smith

{
\seventeenpoint
\centerline{Comparison of Sums of Independent}
\smallskip
\centerline{Identically Distributed Random Variables}
}
\bigskip\bigskip\medskip
\centerline{\caps{S.J.~Montgomery-Smith}%
\footnote{*}%
{Research supported in part by N.S.F.\ Grant
DMS 9201357.}%
}
\smallskip
{
\eightpoint
\centerline{\it Department of Mathematics, University of Missouri,}
\centerline{\it Columbia, MO 65211.}
}
\bigskip\bigskip

\Abstract

Let $S_k$\ be the $k$-th partial sum of Banach space valued
independent identically distributed
random variables.  In this paper, we compare the tail
distribution of $\snormo{S_k}$\ with that of
$\snormo{S_j}$, and 
deduce
some tail distribution maximal inequalities.

\footremark{A.M.S.\ (1991) subject classification: 60-02, 60G50.}

The main result of this paper was inspired by the inequality from \ref[dP--M]
that says that
$\Pr(\snormo{X_1} > t) \le 5 \Pr(\snormo{X_1+X_2} > t/2)$\
whenever $X_1$\ and $X_2$\ are independent identically distributed.  
Such results for $L_p$\ ($p\ge 1$) such as 
$\snormo{X_1}_p \le \snormo{X_1+X_2}_p$\
are straightforward, at least if $X_2$\ has zero expectation.  This inequality
is also obvious if either $X_1$\ is symmetric, or $X_1$\ is real
valued positive.  However, for arbritary random variables, this result
is somewhat surprizing to the author.  Note that
the identically distributed assumption
cannot be dropped, as one could take $X_1 = 1$\ and $X_2 = -1$.

In this paper, we prove a generalization to sums of
arbritarily many independent
identically distributed random variables.  Note that all
results in this paper are true for Banach space valued random variables.
The author would like to thank Victor de la Pe\~na 
for helpful conversations.

\proclaim Theorem 1.  There exist universal constants
$c_1 = 3$\ and $c_2 = 10$\ such that if $X_1,X_2,\dots$\ are independent
identically distributed random variables, and if we set
$$ S_k = \sum_{i=1}^k X_i ,$$
then for $1 \le j \le k$
$$ \Pr(\snormo{S_j} > t) \le c_1 \Pr(\snormo{S_k} > t/c_2) .$$

This result cannot be asymptotically improved.  Consider, for example,
the case where
$X_1 = 1$\ with very small 
probability, and is zero otherwise.  This shows that for the
inequality to be true for all $X_1$, it must be that $c_2$\ must be
larger than some universal constant 
for all $j$\ and $k$.  Also, it is easy to see that $c_1$\ 
must be
larger than some universal constant
because it is easy to select $X_1$\ and $t$\ so that
$\Pr(\snormo{S_j} > t)$\ is close to $1$.

However, Lata{\l}la \ref[L] has been able to obtain the same theorem
with $c_1 = 4$\ and $c_2 = 5$, or $c_1 = 2$\ and $c_2 = 7$.
In the case $j=1$\ and $k=2$, he has shown that
$\Pr(\snormo{X_1} > t) \le 2 \Pr(\snormo{X_1+X_2} > 2t/3)$, and these constants
cannot be improved.

\bigskip

In order to show this result, we will use the following definition.  We 
will say that $x$\ is a {\it $t$-concentration point\/} for a random variable
$X$\ if $\Pr(\snormo{X-x} \le t) > 2/3$.

\proclaim Lemma 2.  If $x$\ is a $t$-concentration point for $X$, and
$y$\ is a $t$-concentration point for $Y$, and $z$\ is a $t$-concentration
point for $X+Y$, then $\snormo{x+y-z} \le 3t$.

\Proof:
$$ \eqalignno{
   \Pr(\snormo{x+y-z} > 3t) &\le
   \Pr(\snormo{X-x + Y-y - (X+Y-z)} > 3t) \cr&\le
   \Pr(\snormo{X-x} > t) + \Pr(\snormo{Y-y} > t) + \Pr(\snormo{X+Y-z} > t) \cr
   &< 1 .\cr}$$
Hence $\Pr(\snormo{x+y-z} \le 3t) > 0$.  
Since $x$, $y$\ and $z$\ are fixed
vectors, the result follows.

\endproof

\proclaim Corollary 3.  If $X_1,X_2,\dots$\ are independent
identically distributed random variables, and if the partial sums
$S_j = \sum_{i=1}^j X_i$\ have $t$-concentration points 
$s_j$\ for $1 \le j \le k$, 
then
$\snormo{k s_j - j s_k} \le 3(k+j)t$.

\Proof:  We prove the result by induction.  It is obvious if
$j=k$.  Otherwise, 
$$ \eqalignno{
   \snormo{ j s_k - k s_j } 
   &\le
   \snormo{ j s_{k-j} - (k-j) s_j } + \snormo{ j s_k - j s_{k-j} - j s_j } \cr
   &\le
   3(k-j + j) t + 3jt = 3(k+j) t .\cr } $$
(The observant reader will notice that we are, in fact, following the
steps of Euclid's algorithm.  The same proof could show
$\snormo{k s_j - j s_k} \le 3 (j+k-2h)t$\ where $h$\ is the highest
common factor of $j$\ and $k$.) 
\endproof

\Proofof Theorem 1:  We consider three cases.
First suppose that $\Pr(\snormo{S_{k-j}} > 9t/10) \le 1/3$.  Note that
$S_k - S_j$\
is independent of $S_j$, and identically distributed
to $S_{k-j}$.  Then
$$ \Pr( \snormo{S_j} > t)
   \le 
   3/2 \Pr( \snormo{S_j} > t \and \snormo{S_k - S_j} \le 9t/10) 
   \le
   3/2 \Pr(\snormo{S_k} > t/10) . $$

For the second case, 
suppose that there is a $1 \le i \le k$\ such that $S_i$\ does
not have any $(t/10)$-concentration point.  Then
$$ \Pr\bigl(\snormo{S_i + X_{i+1}+\dots+X_k} > t/10 \,\big| \,
   \sigma(X_{i+1},\dots,X_k)\bigr) \ge 1/3 ,$$
and hence 
$\Pr(\snormo{S_k} > t/10) \ge 1/3 \ge 1/3 \Pr(\snormo{S_i} > t)$.

Finally, we are left with the third case where 
$\Pr(\snormo{S_{k-j}} > 9t/10) > 1/3$,
and $S_i$\ has a $(t/10)$-concentration point $s_i$\ for all 
$1 \le i \le k$. 
Clearly $\snormo{ s_{k-j} } \ge 8t/10$.  Also, by Corollary~3, 
$$ \snormo{s_k} \ge {k\over k-j} \snormo{s_{k-j}} - {3(2k-j)t\over 10(k-j)}
   \ge {8kt\over 10(k-j)} - {6kt\over 10(k-j)} \ge {2t\over 10} .$$
Therefore, $\Pr(\snormo{S_k} \ge t/10) \ge \Pr(\snormo{S_k-s_k} \le t/10)
\ge 2/3 \ge 2/3 \Pr(\snormo{S_j} > t)$, 
and we are done.
\endproof

One might be emboldened to conjecture the
following.  Suppose that $X_1,X_2,\dots$\ are independent identically
distributed random variables, and that $\alpha_i > 0$.  Let
$$ S_k = \sum_{i=1}^k \alpha_i X_i .$$
Then one might conjecture that there is a universal constant
such that for $1 \le j \le k$
$$\Pr(\snormo{ S_j} > t) \le c\,\Pr(\snormo{S_k} > t/c) .$$

As it turns out, this is not the case.  Let $Y_1,Y_2,\dots$\ be
real valued
independent identically
distributed random variables such that
$$ \eqalignno{
   \Pr( Y_i = N-1) &= 1/N \cr
   \Pr( Y_i = -1 ) &= (N-1)/N .\cr } $$
Then by the central limit theorem, there exists $M \ge N^3$\ such that
$$ \Pr\left( \modo{{1\over M^{2/3}} \sum_{i=1}^M Y_i} > {1\over N}
   \right) \le {1\over N} .$$
Now let $X_i = Y_i + 1/M^{1/3}$, and let
$$ S_M = {1\over M^{2/3}} \sum_{i=1}^M X_i .$$
Then
$ \Pr(\modo{ S_M } > 1/2) \ge 1 - 1/N $,
whereas
$ \Pr(\modo{ S_M + X_{M+1} } > 3/N)  \le 2/N $.

\bigskip

Theorem~1 has several corollaries.

\proclaim Corollary 4.  There is a universal constant $c$\ such that
if $X_1,X_2,\dots$\ are independent
identically distributed random variables, and if we set
$$ S_k = \sum_{i=1}^k X_i ,$$
then
$$ \Pr(\sup_{1 \le j \le k} \snormo{S_j} > t) \le c\, 
   \Pr(\snormo{S_k} > t/c) .$$

Latalo \ref[L] has been able to obtain this result
with $c_1 = 4$\ and $c_2 = 6$, or with $c_1 = 2$\ and $c_2 = 8$.

\Proof:  This follows from Proposition~1.1.1 of \ref[K--W],
that states that if $X_1,X_2,\dots$\ are independent
(not necessarily identically distributed), and if $S_k=\sum_{i=1}^k X_i$,
then
$$ \Pr(\sup_{1 \le j \le k} \snormo{S_j} > t) \le 3
   \sup_{1 \le j \le k}
   \Pr(\snormo{S_j} > t/3) .$$
It is also possible to prove this result directly using the techniques
of the proof of Theorem~1.  The third case only requires that
$\Pr(\normo{S_{k-j}} > 9t/10) > 1/3$\ for one $j = 1$, $2,\dots,$\ $k$.
Hence, for the first case we may assume that $\Pr(\normo{S_k - S_j} > 9t/10) 
\le 1/3$\ for all $1\le j \le k$.  
Let $A_j$\ be the event $\{ \normo{S_i} \le t \hbox{ for all $i < j$\ and }
\normo{S_j} > t\}$.  Then
$$ \Pr(A_j) \le 3/2 \Pr(A_j \and \normo{S_k - S_j} \le 9t/10)
   \le 3/2 \Pr(A_j \and \normo{S_j} > t/10) .$$
Summing over $j$, the result follows.
\endproof

\proclaim Corollary 5.  There is a universal constant $c$\ such that
if $X_1,X_2,\dots$\ are independent
identically distributed random variables, 
and if $\modo{\alpha_i} \le 1$, then
$$ \Pr\left(\Bignormo{ \sum_{i=1}^k \alpha_i X_i } > t\right) \le 
   c\, \Pr\left(\Bignormo{ \sum_{i=1}^k X_i } > t/c\right) .$$

\Proof:  The technique used in this proof is well known
(see for example \ref[KW], Proposition~1.2.1), but is
included for completness.

By taking real and imaginary parts of $\alpha_i$,
we may suppose that the $\alpha_i$\ are real.
Without loss of generality,
$1 \ge \alpha_1 \ge \dots \ge \alpha_k \ge -1$.
Then we may write $\alpha_j = -1 + \sum_{i=j}^k \sigma_i$, where
$\sigma_i \ge 0$.  Thus $\sum_{i=1}^k \modo{\sigma_i} \le 2$, and hence
$$ \eqalignno{
   \Bignormo{ \sum_{j=1}^k \alpha_j X_i }
   &= \Bignormo{ \sum_{j=1}^k \Bigl(-1 + \sum_{i=j}^k \sigma_i\Bigr) 
   X_i } \cr
   &= \Bignormo{ -\Bigl(\sum_{i=1}^k X_i\Bigr)
      + \Bigl(\sum_{j=1}^k \sigma_j \sum_{i=1}^j X_i \Bigr) } \cr
   &\le \Bignormo{ \sum_{i=1}^k X_i }
      + \Bigl( \sum_{j=1}^k \modo{\sigma_j} \Bigr)
        \sup_{1 \le j \le k} \Bignormo{ \sum_{i=1}^j X_i } . \cr }$$
Applying Corollary~4, the result follows.
\endproof

\proclaim Corollary 6.  There are universal constants $c_1$\  
and $c_2$\ such that
if $X_1,X_2,\dots$\ are independent
identically distributed random variables, and if we set
$$ S_k = \sum_{i=1}^k X_i ,$$
then for $1 \le k \le j$
$$ \Pr(\snormo{S_j} > t) \le c_1j/k \, \Pr(\snormo{S_k} > kt/c_2j) .$$

\Proof:  Let $m$\ be the least integer such that $m k \ge j$.  By 
Theorem~1, it follows that $\Pr(\snormo{S_j} > t) 
\le c\,\Pr(\snormo{S_{mk}} > t/c)$.
That $\Pr(\snormo{S_{mk}} > t) \le m \,\Pr(\snormo{S_k} > t/m)$\ is 
straightforward.
\endproof

\bigskip

The example where $X_1$\ is constant shows that $c_2$\ cannot be
made smaller than some universal constant.  The example where
$X_1 = 1$\ with very small
probability and is zero otherwise shows the same is true for $c_1$.

\beginsection References

\references{
dP--M & V.H.~de~la~Pe\~na and S.J.~Montgomery-Smith,\rm\
Bounds on the tail probability of U-statistics and quadratic forms,\sl\
preprint.\cr
K--W & S.~Kwapie\'n and W.A.~Woyczy\'nski,\sl\ 
Random Series and Stochastic
Integrals: Simple and Multiple,\rm\ Birkhauser, NY (1992).\cr
L & R.~Lata{\l}a,\rm\ A paper,\sl\ Warsaw University Preprint (1993).\cr
}

\bye